# For differential equations with $r$ parameters, $2r+1$ experiments are enough for identification

E.D. Sontag[*]

Department of Mathematics
Rutgers University, New Brunswick, NJ 08903
http://www.math.rutgers.edu/~sontag

## Abstract

Given a set of differential equations whose description involves unknown parameters, such as reaction constants in chemical kinetics, and supposing that one may at any time measure the values of some of the variables and possibly apply external inputs to help excite the system, how many experiments are sufficient in order to obtain all the information that is potentially available about the parameters? This paper shows that the best possible answer (assuming exact measurements) is $2r+1$ experiments, where $r$ is the number of parameters.

## 1 Introduction

Suppose that we are given a set of differential equations whose description involves unknown parameters, such as reaction constants in chemical kinetics, resistances in electrical networks, or damping coefficients in mechanical systems. At any time, we may measure the value of some of the variables, or more generally of certain functions of the variables, and we may also apply external inputs to help excite the system so as to elicit more information. Measurements are assumed to be accurate, with no observation noise. We address the following question: *how many experiments are sufficient in order to obtain all the information that is potentially available about the parameters?* The main result is: $2r+1$ *experiments, where $r$ is the number of parameters.*

Questions of this type appear in many areas, and indeed the *identification* (or, when parameters are though of as constant states, the *observation*) problem is one of the central topics in systems and control theory. However, the main motivation for this note arose from recent work on cell signaling pathways. In that field, and in contrast to the standard paradigm in control theory, it is *not* possible to apply arbitrary types of inputs to a system. Often inputs, such as growth factors or hormones, are restricted to be applied as steps of varying durations and amplitudes (or perhaps combinations of a small number of such steps), but seldom does one have the freedom assumed in control-theory studies, where for instance closure of the input class under arbitrary concatenations is needed in order to prove the basic theorems on observability.

The "$2r+1$" expression appears often in geometry and dynamical systems theory. It is the embedding dimension in the "easy" version of Whitney's theorem on representing abstract manifolds as submanifolds of Euclidean space, and it is also the dimension in which $r$-dimensional attractors are embedded, in Takens' classical work. In Aeyels' control theory papers, it is the number of samples needed for observability of generic $r$-dimensional dynamical systems with no inputs. Technically, our problem is quite different, as evidenced by the fact that, in contrast to these studies, which deal with smooth manifolds and systems, we require analytic dependence on parameters, and our main conclusion is false in the more general smooth case. Nonetheless, there are relationships among the topics, which we discuss in the paper.

---

[*]Supported in part by US Air Force Grant F49620-01-1-0063

Let us start to make the problem precise. The system of interest will be assumed to have this form:

$$\dot{z} = f(z, u, x) \qquad (1)$$
$$z(0) = \chi(x) \qquad (2)$$
$$y = h(z, u, x) \qquad (3)$$

where dot indicates derivative with respect to time, and $(dz/dt)(t) = f(z(t), u(t), x)$ depends on a time-invariant parameter $x$ and on the value of the external input $u$ at time $t$. The internal state $z$ of the system, as well as the external inputs, are vector functions, and the parameter $x$ is a constant vector (we prefer to speak of "a vector parameter" as opposed to "a vector of parameters" so that we can later say "two parameters" when referring to two such vectors). The measurements at time $t$ are represented by $y(t) = h(z(t), u(t), x)$; typically, $y$ does not depend directly on inputs nor parameters, and is simply a subset of the state variables $z$, i.e. the function $h$ is a projection. We suppose that the initial state is also parametrized, by a function $\chi$. One more item is added to the specification of the given system: a class of inputs $\mathcal{U}$, meaning a set of functions of time into some set $U$. All definitions will be with regard to time functions $u(t)$ in this set. Before providing more details, let us discuss an example.

## 1.1 An Example

As an illustration, we take the following system:

$$\dot{M} = \frac{E^m}{1 + E^m} - aM, \quad \dot{E} = M - bE \qquad (4)$$

for the production of an enzyme $E$ and the corresponding messenger RNA $M$. The first term for $\dot{M}$ models repression, the positive Hill constant $m$ (not necessarily an integer) specifying the strength of this negative feedback, and the positive constants $a$ and $b$ account for degradation. This is a (nondimensionalized) modification by Griffith [13] of the classical operon model of Goodwin [11]; see also the textbook [6], pp. 208 and 308, as well as the recent paper [19] which describes other variations which are biologically more accurate. In (1), the state $z$ is the vector $(M, E)$. We assume first that there is no true external input to the system, so experiments consist of simply letting the system evolve from its initial state up to certain time $T$, and measuring $M(T)$ at the end of the interval $[0, T]$. (That is, the measured quantity is the amount of RNA; currently gene arrays are used for that purpose.) Experiments differ only by their length $T$. (To fit in the general framework, where inputs are allowed, we may simply take the set $U$ of input values to be, for instance, $\{0\}$, so that the set $\mathcal{U}$ consists of just one function, the input $u \equiv 0$.) We take as parameters for this problem the possible 5-dimensional positive vectors:

$$x = (M_0, E_0, m, a, b)$$

which list the initial conditions as well as the three constants. Thus $\chi(x)$ is just the vector $(M_0, E_0)$, and $f(M, E, M_0, E_0, m, a, b) = (\frac{E^m}{1+E^m} - aM, M - bE)$. As output function $h$, we take $h(M, E) = M$. It is not possible to identify all parameters in this example: for any positive numbers $k$ and $\ell$, the parameter

$$x = \left(1, k, \ell, \frac{k^\ell}{1 + k^\ell}, \frac{1}{k}\right)$$

gives rise to the same output $M(t) \equiv 1$. Our general theorem will imply that a set of measurements taken at a random choice of $2r+1 = 11$ instants is sufficient in order to distinguish between any two parameters which give rise to different output functions $M(t)$. (Less than $2r+1$ measurements may be enough in any given example; the $2r+1$ bound is a very general upper bound, which is best possible in the sense that for *some* systems, no less will do.)

More interestingly, let us now suppose that in our experiments we can affect the degradation rate $b$. Specifically, suppose that another substance, which binds to (and hence neutralizes) the enzyme $E$, is added at a concentration $u$ which we can choose, and mass action kinetics controls the binding. We pick units of $u$ such that the new equations become

$$\dot{M} = \frac{E^m}{1 + E^m} - aM, \quad \dot{E} = M - bE - uE. \tag{5}$$

The parameters are the same as before, and neither $\chi$ or $h$ need to be modified; the only change is that now $f(M, E, M_0, E_0, m, a, b, u) = (\frac{E^m}{1+E^m} - aM, M - bE - uE)$ because there is an explicit input in the system description. We will assume that the concentration $u$ is kept constant at a value $u_0$ for the duration of the experiment, so that the set of possible experiments is now specified by two numbers instead of one: a pair $(u_0, T)$ consisting of the concentration $u_0 \geq 0$ (so, $u_0 = 0$ means no neutralization) and the length of the time interval $T$ at the end of which we measure $M(T)$. Here the set $\mathcal{U}$ consists of all possible nonnegative constant mappings; $[0, \infty) \to \mathbb{R}$; thus $(u_0, T)$ represents the input function $u = u_0$ used on the interval $[0, T]$. In this case, it turns out, every pair of distinct parameters is distinguishable (see Section 5 for the calculation), so the general theorem implies that a set of measurements taken at a random choice of $2r+1 = 11$ instants is sufficient in order to identify $(M_0, E_0, m, a, b)$. As yet another variation of this example, one might be able, for instance, to change the concentration $u(t)$ in a linearly increasing fashion: $u(t) = u_0 + u_1 t$, with $u_0 \geq 0$ and $u_1 \geq 0$. Now the experiments are specified by triples $(u_0, u_1, T)$, where again $T$ stands for the time at which we take the measurement, and $\mathcal{U}$ consists of all linear functions as above.

## 1.2 The Main Result

We will define an "analytically parametrized system" as one in which the system of differential equations, the initial state map $\chi$, and the observation map $h$ are all expressed as real-analytic functions of states, parameters, and inputs, and also inputs depend analytically on a finite number of parameters. Recall that real-analytic maps are those that may be represented locally about each point of their domain by a convergent power series. This includes most reasonable nonlinearities one may think of: rational functions, roots and exponents (as long as away from singularities), trigonometric functions, logarithms and exponentials, and so forth, but not switching discontinuities nor smooth "patchings" between discontinuous pieces. (The results are definitely not true if analyticity is relaxed to just infinite differentiability, as we shall show by counterexample.)

We will also define precisely what we mean by "a randomly chosen set of $2r+1$ experiments" and by "distinguishability" of parameters. All these definitions express the intuitive ideas conveyed by the terms, but are necessarily somewhat technical so we defer them until after the statement of our main theorem.

**Theorem 1** *Assume given an analytically parametrized system, and let $r$ be the dimension of its parameter space. Then, for any randomly chosen set of $2r+1$ experiments, the following property holds: for any two parameters that are distinguishable, one of the experiments in this set will distinguish them.*

Let us now make the terms precise. In order to properly define analytic maps, we need that states $z(t)$, parameters $x$, input values $u(t)$, and measurement values $y(t)$ all belong to analytic manifolds. Examples of such manifolds are all open subsets of an Euclidean space, which are indeed the most common situation in applications; for instance, in biological applications states and parameters are usually given by vectors with strictly positive entries. But using manifolds more general than open subsets of Euclidean spaces allows one to model constraints (such as the fact that a state may be an angle, i.e. an element of the unit circle), and adds no complexity to the proof of the main result. We do not loose generality, however, in assuming that $y(t) \in \mathbb{R}^p$ for some integer $p$, since every analytic manifold can be embedded in some Euclidean space.

Formally, we define *analytically parametrized systems* as 9-tuples:

$$\Sigma = (M, U, \mathbb{X}, \Lambda, p, f, \chi, h, \mu)$$

where

- the *state-space* $M$, *input-value space* $U$, *parameter space* $\mathbb{X}$, and *experiment space* $\Lambda$ are real-analytic manifolds,

- $p$ is a positive integer

- $f : M \times U \times \mathbb{X} \to TM$ (the tangent bundle of $M$), $\chi : \mathbb{X} \to M$, $h : M \times U \times \mathbb{X} \to \mathbb{R}^p$, and $\mu : \Lambda \times \mathbb{R} \to U$ are real-analytic maps, and

- $f$ is a vector field on $M$, that is, $f(z, u, x) \in T_z M$ for each $(z, u, x)$.

An additional technical assumption, completeness, will be made below.

Given a system $\Sigma$, we associate to it its *response*, the mapping

$$\beta_\Sigma : \mathbb{X} \times \Lambda \to \mathbb{R}^p$$

obtained, for each system parameter $x \in \mathbb{X}$ and experiment $\lambda \in \Lambda$, by solving the initial-value problem (1)-(2) with the input $u(t) = \mu(\lambda, t)$, and then evaluating the output (3) at the final time $t = 1$ (We are arbitrarily normalizing the time interval to $[0, 1]$, but varying lengths can be easily incorporated into the formalism, as we discuss later.) Thus, the class of inputs $\mathcal{U}$ is the set of all maps $t \mapsto \mu(\lambda, t)$. More precisely: we consider the solution $z(\cdot)$ of the differential equation $(dz/dt)(t) = g(t, z(t))$ with initial condition $z(0) = \chi(x)$, where $g$ is the time-varying vector field defined by $g(t, \zeta) = f(\zeta, \mu(\lambda, t), x)$. Since $f(\zeta, \mu(\lambda, t), x)$ depends analytically on $t$, $\zeta$, $\lambda$, and $x$, such a solution exists at least for all $t \approx 0$, and it depends analytically on $\lambda$ and $x$. (This is a standard fact about differential equations; see e.g. Proposition C.3.12 in [25], viewing parameters as constant additional states, and note that the manifold case follows easily from the Euclidean case.) We make the following *completeness* assumption: the solution $z$ exists on the interval $[0, 1]$. Now we define $\beta_\Sigma(x, \lambda) := h(z(1), \mu(\lambda, 1), x)$. Thus, $\beta_\Sigma$ is a real-analytic function.

We define two parameters $x_1$ and $x_2$ to be *distinguishable* if there exists some experiment $\lambda$ which distinguishes between them:

$$\beta_\Sigma(x_1, \lambda) \neq \beta_\Sigma(x_2, \lambda).$$

We next discuss the notion of "random" set of $2r+1$ experiments. In general, we will say that a property holds *generically* (often the term "residual" is used for this concept) on a topological space $S$ if the set of elements $S_0 \subseteq S$ for which the property holds contains the intersection of a countable family of open dense sets. Such sets are "large" in a topological sense, and for all the spaces that we consider, the "Baire property" holds: generic subsets $S_0$ are dense in $S$. When $S$ is a manifold, there is a well-defined concept of measure zero subset, see e.g. [10, 14]; in that case we will say that a subset $S_0 \subseteq S$ has *full measure* if its complement $S \setminus S_0$ has measure zero. (See e.g. [17] for a comparison of these two alternative concepts of "large" subset of $S$. For an extension to infinite dimensional spaces of the notion of full measure, called "prevalence," see [15].) We will show that sets of $2r+1$ experiments which distinguish are *both* generic *and* of full measure.

In general, for any set $M$ and any positive integer $s$, as in [10] we denote by $M^{(s)}$ the subset of $M^s$ made up of all sequences $(\xi_1, \ldots, \xi_s)$ consisting of distinct elements ($\xi_i \neq \xi_j$ for each $i \neq j$). Then, we say that a property (P) of $q$-element sets of experiments holds for a random set of $q$ experiments, where $q$ is a positive integer, provided that

$$\mathcal{G}_q = \left\{ (\lambda_1, \lambda_2, \ldots, \lambda_q) \mid \text{(P) holds for the set } \{\lambda_1, \lambda_2, \ldots, \lambda_q\} \right\}$$

is generic and of full measure in $\Lambda^{(q)}$. Now all the terms in the statement of Theorem 1 have been defined.

**Remark 1.1** In applications, often we may not want to restrict the space $\Lambda$ which parametrizes experiments to have to be a manifold. For instance, in the examples in Section 1.1, we took constant or linear nonnegative inputs, meaning that $\Lambda = [0, \infty)$ or $[0, \infty)^2$ respectively. Such generalizations present no difficulty, and can be handled by the theory in several alternative ways. The simplest is just to pick a different parametrization of inputs. For example, we may write constant nonnegative inputs as $u_0^2$, with $u_0 \in \Lambda = \mathbb{R}$ and linear nonnegative inputs as $u_0^2 + u_1^2 t$, with $(u_0, u_1) \in \Lambda = \mathbb{R}^2$. Another possibility is to use a smaller $\Lambda$ which is dense in the original one; for instance, we may use $\Lambda = (0, \infty)$ for the constant input case, and note that, by continuity of $\beta_\Sigma(x, \lambda)$ on $\lambda$, distinguishability is not affected when using such restricted experiments. Finally, the best alternative would be to prove the theorem for more general sets $\Lambda$ (and $\mathbb{X}$ too), namely arbitrary "subanalytic" subsets of analytic manifolds; we did not do so in order to keep the presentation as simple as possible. □

**Remark 1.2** We normalized the time interval to $[0, 1]$. However, we can easily include the case in which observations may be performed on the system at different times. We simply view measurements taken at different times as different experiments, adding the final time $T$ as a coordinate to the specification of experiments, just as we did in the examples discussed in Section 1.1. Formally, we define the map

$$\beta_\Sigma^* \,:\, \mathbb{X} \times \Lambda \times (0,\infty) \to \mathbb{R}^p \,:\, (x,\lambda,T) \mapsto h(z(T), \mu(\lambda,T), x)$$

and say that two parameters $x_1$ and $x_2$ are *distinguishable on varying intervals* if there exist $\lambda \in \Lambda$ and $T > 0$ such that $\beta_\Sigma^*(x_1,\lambda,T) \neq \beta_\Sigma^*(x_2,\lambda,T)$. Theorem 1 gives as a corollary that, for a random set of pairs of the form $\{(\lambda_1,T_1),(\lambda_2,T_2),\ldots,(\lambda_{2r+1},T_{2r+1})\}$, if any two parameters $x_1$ and $x_2$ are distinguishable on varying intervals, then $\beta_\Sigma^*(x_1,\lambda_i,T_i) \neq \beta_\Sigma^*(x_2,\lambda_i,T_i)$ for at least one of the $2r+1$ pairs in this set. To prove this corollary, we introduce a new set of experiments $\widetilde{\Lambda} := \Lambda \times (0,\infty)$, and a new analytically parametrized system

$$\widetilde{\Sigma} = \left( M, \widetilde{U}, \mathbb{X}, \widetilde{\Lambda}, p, \widetilde{f}, \chi, \widetilde{h}, \widetilde{\mu} \right)$$

with the property that $\beta_\Sigma^*(x,\lambda,T) = \beta_{\widetilde{\Sigma}}(x,(\lambda,T))$ for all $(x,\lambda,T)$. Clearly, the desired conclusion follows for $\Sigma$ when the theorem is applied to $\widetilde{\Sigma}$. To define $\widetilde{\Sigma}$, we note that we must have

$$h(z(T), \mu(\lambda,T), x) \;=\; \widetilde{h}(\widetilde{z}(1), \widetilde{\mu}((\lambda,T),1), x) \tag{6}$$

for all $(x,\lambda,T)$, where $\widetilde{z}$ solves $(d\widetilde{z}/dt)(t) = \widetilde{f}(\widetilde{z}, \widetilde{\mu}((\lambda,T),t), x)$ with initial condition $\widetilde{z}(0) = \chi(x)$. This can be accomplished by reparametrizing time according to the experiment duration:

$$\widetilde{U} := (0,\infty) \times U, \quad \widetilde{\mu}((\lambda,T),t) := (T, \mu(\lambda,Tt)),$$

$$\widetilde{f}(z,(u_0,u),x) := u_0 f(z,u,x), \quad \widetilde{h}(z,(u_0,u),x) := h(z,u,x)\,.$$

It is easy to verify that (6) holds with these choices. □

The rest of this note is organized as follows. Section 2 studies the more abstract problem of distinguishability for "response" maps $\beta$ which do not necessarily arise from systems $\Sigma$, and presents a general $2r+1$ theorem for such responses. That result already appeared, expressed in a slightly different manner, as a technical step in [23], so we do not include all the details of the proof, and in particular the technical material on analytic functions, which can be found in that reference. Also in Section 2, we present results for merely smooth mappings: one showing that $r$ experiments are enough for local identification, and another one on genericity. Section 3 specializes to responses linear in inputs, a class of responses which is of interest because of relations to Whitney's embedding theorem and data reduction (cf. Section 6.6), and, in Section 3.1, we provide a nontrivial example in this class showing that the bound $2r+1$ cannot be improved in the analytic case. Section 4 shows how the results apply to the system case, proving Theorem 1. Section 5 completes the discussion of the second example in Section 1.1, showing that all pairs of parameters are distinguishable. Finally, we close with Section 6, where we provide general comments and discuss relations to other work.

## 2 Parameter Distinguishability for Maps

In this section, we consider maps, which we call *responses*

$$\beta \,:\, \mathbb{X} \times \mathbb{U} \to \mathbb{R}$$

where $\mathbb{X}$ and $\mathbb{U}$ are two differentiable (connected, second countable) manifolds. We may view a map $\beta(x,\cdot) : \mathbb{U} \to \mathbb{R}$, for each fixed (vector) "parameter" $x$, as the (scalar) response of a system to "inputs" in $\mathbb{U}$ (in the application to systems, these will be elements in $\Lambda$ which parametrize continuous-time inputs). In typical applications, $\mathbb{X}$ is an open subset of an Euclidean space $\mathbb{R}^r$ and $\mathbb{U}$ is an open subset of some $\mathbb{R}^m$. In general, we call $\mathbb{X}$ the *parameter space* and $\mathbb{U}$ the *input space*, and use respectively $r$ and $m$ to denote their dimensions.

The results of most interest here will be those in which $\beta$ is a (real-)analytic map (that is, $\beta$ may be represented by a locally convergent power series around each point in $\mathbb{X} \times \mathbb{U}$), in which case we assume implicitly that $\mathbb{X}$ and $\mathbb{U}$ are analytic manifolds. However, we also will make some remarks which apply to the more general cases when $\beta$ is merely smooth, i.e. infinitely differentiable (and $\mathbb{X}$ and $\mathbb{U}$ are smooth manifolds), or even just continuously differentiable.

Given a response $\beta$ and any subset of inputs $\mathbb{U}_0 \subseteq \mathbb{U}$, two parameters $x_1$ and $x_2$ are said to be *indistinguishable by inputs in* $\mathbb{U}_0$, and we write

$$x_1 \underset{\mathbb{U}_0}{\sim} x_2 \,,$$

if

$$\beta(x_1, u) \;=\; \beta(x_2, u) \quad \forall\, u \in \mathbb{U}_0\,.$$

If this property holds with $\mathbb{U}_0 = \mathbb{U}$, we write $x_1 \sim x_2$ and simply say that $x_1$ and $x_2$ are *indistinguishable*; this means that $\beta(x_1, u) = \beta(x_1, u)$ for all $u \in \mathbb{U}$: the response is the same, for all possible inputs, whether the parameter is $x_1$ or $x_2$. If, instead, there exists some $u \in \mathbb{U}$ such that $\beta(x_1, u) \neq \beta(x_1, u)$, we say that $x_1$ and $x_2$ are *distinguishable*, and write $x_1 \not\sim x_2$.

We say that a set $\mathbb{U}_0$ is a *universal distinguishing set* if

$$x_1 \underset{\mathbb{U}_0}{\sim} x_2 \quad \Longleftrightarrow \quad x_1 \sim x_2$$

which means that if two parameters can be distinguished at all, then they can be distinguished on the basis of inputs taken from the subset $\mathbb{U}_0$ alone.

A useful notation is as follows. For each fixed positive integer $q$, we extend the function $\beta : \mathbb{X} \times \mathbb{U} \to \mathbb{R}$ to a function

$$\beta_q \;:\; \mathbb{X} \times \mathbb{U}^q \to \mathbb{R}^q \;:\; (x, u_1, \ldots, u_q) \mapsto \begin{pmatrix} \beta(x, u_1) \\ \vdots \\ \beta(x, u_q) \end{pmatrix}$$

and, with some abuse of notation, we drop the subscript $q$ when clear from the context. Then, saying that the finite subset $\mathbb{U}_0 = \{u_1, \ldots, u_q\}$ is a universal distinguishing set amounts to the following property holding for all $x_1, x_2 \in \mathbb{X}$:

$$(\exists\, u \in \mathbb{U})\, (\beta(x_1, u) \neq \beta(x_2, u)) \quad \Rightarrow \quad \beta(x_1, u_1, \ldots, u_q) \neq \beta(x_2, u_1, \ldots, u_q)\,. \tag{7}$$

## 2.1 Global Analytic Case: $2r+1$ Experiments are Enough

We now turn to the main theorem of this section. It states that $2r+1$ (recall $r = \dim \mathbb{X}$) experiments are sufficient for distinguishability, and, moreover, that a random set of $2r+1$ inputs is good enough, if $\beta$ is analytic. Later, we show that the bound $2r+1$ is best possible. This theorem was already proved, using slightly different terminology, in [23]. We provide here the outline of the proof, but for the main technical step, concerning real-analytic manifolds, we will refer the reader to [23].

We express our result in terms of the set $\mathcal{G}_{\beta,q}$ consisting of those ("good") sequences $(u_1, \ldots, u_q) \in \mathbb{U}^{(q)}$ which give rise, for a given response $\beta$, to universal distinguishing sets $\mathbb{U}_0 = \{u_1, \ldots, u_q\}$ of cardinality $q$, and its complement $\mathbb{U}^{(q)} \setminus \mathcal{G}_{\beta,q}$, the ("bad") sequences:

$$\mathcal{B}_{\beta,q} := \left\{ w \in \mathbb{U}^{(q)} \mid \exists\, x_1, x_2 \in \mathbb{X} \text{ s.t. } x_1 \not\sim x_2 \text{ and } \beta(x_1, w) = \beta(x_2, w) \right\}. \tag{8}$$

**Theorem 2** *Assume that $\beta$ is analytic. Then the set $\mathcal{B}_{\beta, 2r+1}$ is a countable union of embedded analytic submanifolds of $\mathbb{U}^{2r+1}$ of positive codimension. In particular, $\mathcal{G}_{\beta, 2r+1}$ is generic and of full measure in $\mathbb{U}^{2r+1}$.*

We give the proof after establishing a general technical fact. For each fixed positive integer $q$, we introduce the following set:

$$\mathcal{P}_{\beta,q} := \left\{ \left((x_1, x_2), w\right) \in \mathbb{X}^{(2)} \times \mathbb{U}^{(q)} \mid x_1 \not\sim x_2 \text{ and } \beta(x_1, w) = \beta(x_2, w) \right\}. \tag{9}$$

This is a "thin" (possibly empty) subset of $\mathbb{X}^{(2)} \times \mathbb{U}^{(q)}$:

**Lemma 2.1** The set $\mathcal{P}_{\beta,q}$ is a countable union of submanifolds of $\mathbb{X}^{(2)} \times \mathbb{U}^{(q)}$ each of which has dimension $\leq qm + 2r - q$.

*Proof.* The proof is included in [23]; let us recall the main steps. We first consider the set of pairs of parameters which can be distinguished from each other:

$$\mathcal{X} := \{(x_1, x_2) \in \mathbb{X}^{(2)} \mid x_1 \not\sim x_2\}. \tag{10}$$

If this set is empty then $\mathcal{P}_{\beta,q}$ is also empty, and we are done. Otherwise, $\mathcal{X}$ is an open subset of $\mathbb{X}^{(2)}$, and hence an open subset of $\mathbb{X}^2$, and is thus a manifold of dimension $2r$. We let $\mathcal{U}(x_1, x_2)$ be the set consisting of those $u \in \mathbb{U}$ such that $x_1 \underset{u}{\sim} x_2$. If $(x_1, x_2) \in \mathcal{X}$, then $\mathcal{U}(x_1, x_2)$ is an analytic subset (a set defined by zeroes of analytic functions) of $\mathbb{U}$ of dimension at most $m - 1$, since it is the set where the nonzero analytic function $\beta(x_1, u) - \beta(x_2, u)$ vanishes, and $\mathbb{U}$ is connected. Therefore, its Cartesian product $(\mathcal{U}(x_1, x_2))^q$ is an analytic subset of $\mathbb{U}^q$ of dimension at most $q(m-1)$ (Proposition A.2, Part 3, in [23]). Next, for each $(x_1, x_2) \in \mathcal{X}$, we consider the following subset of $\mathbb{U}^q$:

$$\mathcal{T}(x_1, x_2) = \{w \in \mathbb{U}^{(q)} \mid \beta(x_1, w) = \beta(x_2, w)\}.$$

Clearly, $\mathcal{T}(x_1, x_2)$ is a semianalytic subset, i.e. a set defined by analytic equalities (responses are equal) and inequalities (the coordinates of $w$ are distinct). The key point is that $\mathcal{T}(x_1, x_2)$ has dimension at most $q(m-1)$, because it is a subset of $(\mathcal{U}(x_1, x_2))^q$.

The set $\mathcal{P}_{\beta,q}$ is also semianalytic. Let $\pi_1 : \mathbb{X}^{(2)} \times \mathbb{U}^{(q)} \to \mathbb{X}^{(2)}$ be the projection on the first factor. For each $(x_1, x_2) \in \mathcal{X}$, $\pi_1^{-1}(x_1, x_2) \bigcap \mathcal{P}_{\beta,q} = \mathcal{T}(x_1, x_2)$ has dimension at most $q(m-1)$. Applying then Proposition A.2, Part 2, in [23], and using that $\mathcal{X}$ is an analytic manifold of dimension $2r$, it follows that $\dim \mathcal{P}_{\beta,q} \leq 2r + q(m-1) = 2r + qm - q$. It is known from stratification theory (cf. [3, 8, 27], and the summary in the Appendix to [23]) that any semianalytic set is a countable union of embedded analytic submanifolds, so the Lemma is proved. ∎

**Proof of Theorem 2**

Letting $\pi_2 : \mathbb{X}^{(2)} \times \mathbb{U}^{(q)} \to \mathbb{U}^{(q)}$ be the projection on the second factor, we have that $\pi_2(\mathcal{P}_{\beta,q}) = \mathcal{B}_{\beta,q}$. Again from stratification theory, we know that the image under an analytic map of a countable union of (embedded, analytic) submanifolds of dimension $\leq p$ is again a countable union of submanifolds of dimension $\leq p$. In the particular case that $q = 2r+1$, and applying Lemma 2.1, this means that $\mathcal{B}_{\beta,2r+1}$ is a countable union of submanifolds $M_i$ of $\mathbb{U}^{2r+1}$ of dimension $\leq p$, where $p \equiv (2r+1)m + 2r - (2r+1) = (2r+1)m - 1$. Since $\dim \mathbb{U}^{2r+1} = (2r+1)m$, each $M_i$ has positive codimension. Embedded submanifolds are, in local coordinates, proper linear subspaces (and one may cover them by countably many charts), so the complements of each $M_i$ are generic and of full measure, from which it follows that $\mathcal{G}_{\beta,2r+1}$ is generic and of full measure, as wanted. ∎

## 2.2 Local Case: $r$ Experiments are Enough

If only distinguishability of parameters *near* a given parameter is desired, then $r$, rather than $2r+1$, experiments are sufficient, at least for a generic subset of parameters. To make this fact precise, we introduce local versions of the sets of "good" and "bad" inputs. For each open subset $V \subseteq \mathbb{X}$, we let

$$\mathcal{B}_{\beta,q,V} := \left\{ w \in \mathbb{U}^{(q)} \mid \exists x_1, x_2 \in V \text{ s.t. } x_1 \not\sim x_2 \text{ and } \beta(x_1, w) = \beta(x_2, w) \right\}$$

and $\mathcal{G}_{\beta,q,V} = \mathbb{U}^{(q)} \setminus \mathcal{B}_{\beta,q,V}$. When $V = \mathbb{X}$, these are $\mathcal{B}_{\beta,q}$ and $\mathcal{G}_{\beta,q}$ respectively.

**Proposition 2.2** Assume that $\beta(\cdot, u)$ is continuously differentiable for each $u \in \mathbb{U}$. Then there exists an open dense subset $\mathbb{X}_0$ of $\mathbb{X}$ with the following property: for each $x \in \mathbb{X}_0$ there is some neighborhood $V$ of $x$ such that $\mathcal{G}_{\beta,r,V}$ is nonempty.

To show the proposition, we first introduce the notion of a nonsingular parameter. For each $x \in \mathbb{X}$, each positive integer $q$, and each $w \in \mathbb{U}^{(q)}$, we define $\rho(x, w)$ to be the rank of the differential of $\beta_q(\xi, w)$ with respect to $\xi$, evaluated at $\xi = x$, and for each $x \in \mathbb{X}$ define

$$\rho(x) := \max \left\{ \rho(x, w) \mid w \in \mathbb{U}^{(q)},\ q \geq 1 \right\}.$$

Observe that the maximum is always achieved at some $w \in \mathbb{U}^{(\rho(x))}$, since if $\rho(x, u_1, \ldots, u_p) = q$ then there must exist some $q$-element subset of $\{u_1, \ldots, u_p\}$ for which the rank is also $q$. We will say that a parameter value $x$ is *nonsingular* provided that this maximal rank is locally constant at $x$, that is, there is some neighborhood $V$ of $x$ in $\mathbb{X}$ such that $\rho(\xi) = \rho(x)$ for every $\xi \in V$. *The set $\mathbb{X}_{\text{NS}}$ of nonsingular parameter values is an open and dense subset of $\mathbb{X}$.* Openness is clear from the definition and density follows from this argument: suppose there would exist an open subset $W \subseteq \mathbb{X}$ with $W \cap \mathbb{X}_{\text{NS}} = \emptyset$; now pick a point $x \in W$ at which $\rho(x)$ is maximal with respect to points in $W$; since $\beta$ is continuously differentiable on $x$, there is a neighborhood of $x$ in $W$ where the rank is at least equal to $\rho(x)$, and thus, by maximality, it is equal to $\rho(x)$; so $x \in \mathbb{X}_{\text{NS}}$, a contradiction. Proposition 2.2 is then a consequence of the following, using $\mathbb{X}_0 = \mathbb{X}_{\text{NS}}$, and because $\rho(x) \leq r$ for all $x$:

**Lemma 2.3** Pick any nonsingular parameter value $x$. Then there is some open neighborhood $V$ of $x$ such that $\mathcal{G}_{\beta,\rho(x),V} \neq \emptyset$.

*Proof.* We must prove that there is an open neighborhood $V$ of $x$ and a subset $\mathbb{U}_0$ of $\mathbb{U}$ of cardinality $q = \rho(x)$ such that $\mathbb{U}_0$ is a universal distinguishing set for parameters in $V$, that is,

$$x_1 \underset{\mathbb{U}_0}{\sim} x_1 \text{ and } x_1, x_2 \in V \Rightarrow x_1 \sim x_2 \, .$$

If $w = (u_1, \ldots, u_q) \in \mathbb{U}^{(q)}$ is such that $q = \rho(x) = \rho(x, w)$, a routine functional dependence argument implies that $\mathbb{U}_0 = \{u_1, \ldots, u_q\}$ has the desired properties, as follows. We pick a neighborhood $V$ of $x$ so that, for any other $u \in \mathbb{U}$, the rank of $\beta_{q+1}(\xi, w, u)$ with respect to $\xi$ is also constantly equal to $q$ on $V$ (nonsingularity of $x$). Consider the mapping

$$K \,:\, V \to \mathbb{R}^q \,:\, \xi \mapsto \beta_q(\xi, w) \, .$$

Since the rank of the differential of $K$ is constant on $V$, by the Rank Theorem and after shrinking $V$ if needed, we know that there exist two diffeomorphisms $S : \mathbb{R}^r \to V$ and $T : W \to \mathbb{R}^q$, where $W$ is the image $K(V)$, such that $T \circ K \circ S : \mathbb{R}^r \to \mathbb{R}^q$ is the canonical projection $(z_1, \ldots, z_r) \to (z_1, \ldots, z_q)$ on the first $q$ coordinates.

Now pick any other $u \in \mathbb{U}$, introduce

$$H \,:\, V \to \mathbb{R}^{q+1} \,:\, \xi \mapsto \beta_{q+1}(\xi, u_1, \ldots, u_q, u) = \begin{pmatrix} \beta_q(\xi, w) \\ \beta(\xi, u) \end{pmatrix} \, ,$$

and let $F : \mathbb{R}^r \to \mathbb{R}^{q+1}$ be the map obtained as the composition $(T \times I) \circ H \circ S$, where $T \times I$ maps $(a, b) \mapsto (T(a), b)$. Since $S : \mathbb{R}^r \to V$ and $(T \times I) : W \to \mathbb{R}^q$ are diffeomorphisms and $H$ has constant rank, also $F$ has constant rank $q$. So, since the Jacobian of $F$ has the block form $\begin{pmatrix} I & 0 \\ * & M \end{pmatrix}$, where $*$ is irrelevant and $M$ is the Jacobian of $\beta(S(\cdot), u)$ with respect to the variables $z_{q+1}, \ldots, z_r$, it follows that $M \equiv 0$. In other words, $\beta(S(\cdot), u)$ must be independent of these variables, so there exists a $\varphi : \mathbb{R}^q \to \mathbb{R}$ such that $\beta(S(z), u) = \varphi(z_1, \ldots, z_q)$ for all $z = (z_1, \ldots, z_r) \in \mathbb{R}^r$. Since $T \circ K \circ S$ is the projection on the first $q$ variables, this means that $\beta(S(z), u) = \varphi((T \circ K \circ S)(z)) = \varphi(T(\beta(S(z), w)))$ for all $z \in \mathbb{R}^r$. As $S$ is onto, this is equivalent to:

$$\beta(\xi, u) \;=\; \varphi\left(T(\beta(\xi, w))\right)$$

for every $\xi \in V$. Now let us suppose that $x_1 \underset{\mathbb{U}_0}{\sim} x_2$. By definition, $\beta(x_1, w) = \beta(x_2, w)$, and this in turn implies

$$\beta(x_1, u) = \varphi(T(\beta(x_1, w))) = \varphi(T(\beta(x_2, w))) = \beta(x_2, u) \, .$$

As $u$ was arbitrary, we conclude $x_1 \sim x_2$. ∎

**Remark 2.4** In general, $\mathbb{X}_0$ must be a proper subset of $\mathbb{X}$. A counterexample to equality would be provided by $\beta(x, u) = f(x) \cdot u$ (dot indicates inner product, see Section 3 for responses of this special form), where $x \in \mathbb{R}$ and $f : \mathbb{R} \to \mathbb{R}^2$ parametrizes a figure-8 curve. Around the parameter value that corresponds to the crossing, no single $u$ will distinguish. (We omit details.) □

## 2.3 Smooth Case: $2r+1$ Experiments are Enough, for Generic $\beta$

For infinitely differentiable but non-analytic $\beta$, the result on existence of universal distinguishing subsets of cardinality $2r+1$ does not generalize. As a matter of fact, one can exhibit a $\beta$ of class $\mathcal{C}^\infty$, with $r = 1$, with the property that there is *no finite* universal distinguishing set: take $\mathbb{X} = \mathbb{U} = (0, \infty)$ and $\beta(x, u) := \gamma(x - u)$, where $\gamma$ is a smooth map which is nonzero on $(-\infty, 0)$ and zero elsewhere. For this example, every two parameters are distinguishable (if $x \neq y$ then picking $u := (x + y)/2$ results in $\beta(x, u) \neq 0 = \beta(y, u)$). To see that there is no finite universal distinguishing set, suppose that $\mathbb{U}_0 \subseteq \mathbb{U}$ is such a set; then picking $x, y \geq \max\{u \mid u \in \mathbb{U}_0\}$ we have that $\beta(x, u) = 0 = \beta(y, u)$ for all $u \in \mathbb{U}_0$, a contradiction. Observe that lack of compactness is not the reason that finite distinguishing subsets fail to exist in the non-analytic case. Let us sketch an example with $\mathbb{X} = \mathbb{S}^1$ and $\mathbb{U} = (-\pi, \pi)$. For each $\varphi \in \mathbb{U}$, we take $\beta(\cdot, \varphi)$ to be a "bump function" centered at $e^{i\varphi}$, with value $= 1$ there and less than one elsewhere, and supported on a compact set which does not contain $x = 1$. (Such a choice can be made smoothly on $x$ and $u$ simultaneously.) If $x \neq y$ are elements of $\mathbb{S}^1$, we can always distinguish them by an appropriate $u$ (if $x \neq 1$, we may take $\varphi$ to be the argument of $x$, so that $\beta(x, u) = 1$ and $\beta(y, u) < 1$, and if $x = 1$ and $y \neq 1$, we may take $\varphi$ to be the argument of $y$). However, given any finite set of $u$'s, the union of the supports of the corresponding bump functions is still a compact subset of $\mathbb{S}^1$ which does not include 1, so there is some $x \neq 1$ such that $\beta(x, u) = 0$ for all $u$ in this set, and hence $x$ is not distinguishable from 1.

It is possible, however, to provide a result that holds *generically* on $\mathcal{C}^\infty(\mathbb{X} \times \mathbb{U}, \mathbb{R})$, the set of smooth (that is, infinitely differentiable) responses $\beta : \mathbb{X} \times \mathbb{U} \to \mathbb{R}$, endowed with the Whitney topology. (Density in the Whitney topology means approximability in a very strong sense; see e.g. [10, 14] for details.) Generic results in this sense are not terribly interesting, and are fairly meaningless in applications, since a response "close" to a given one may not have any interesting structure, but we present the result nonetheless for completeness. The analogue of Theorem 2 is as follows.

**Proposition 2.5** For generic $\beta \in \mathcal{C}^\infty(\mathbb{X} \times \mathbb{U}, \mathbb{R})$, $\mathcal{G}_{\beta, 2r+1}$ is generic and of full measure in $\mathbb{U}^{2r+1}$.

We need this analog of Lemma 2.1:

**Lemma 2.6** For generic $\beta \in \mathcal{C}^\infty(\mathbb{X} \times \mathbb{U}, \mathbb{R})$, $\mathcal{P}_{\beta, q}$ is a submanifold of $\mathbb{X}^{(2)} \times \mathbb{U}^{(q)}$ of dimension $\leq qm + 2r - q$.

Let us show first how Proposition 2.5 follows from here. In the particular case $q = 2r+1$, the Lemma gives that $\dim \mathcal{P}_{\beta, q} \leq (2r+1)m - 1$. The projection $\pi_2 : \mathbb{X}^{(2)} \times \mathbb{U}^{(q)} \to \mathbb{U}^{(q)}$ restricts to a smooth map $f : \mathcal{P}_{\beta, q} \to \mathbb{U}^{(q)}$ with image $f(\mathcal{P}_{\beta, q}) = \mathcal{B}_{\beta, q}$. In general, the Morse-Sard Theorem (as stated e.g. in [14], Theorem 3.1.3) says that if $f : M \to N$ is smooth then $N \setminus f(\Sigma_f)$ is generic and has full measure, where $\Sigma_f$ is the set of critical points of $f$ (differential is not onto). In our case, $\dim \mathcal{P}_{\beta, q} < (2r+1)m = \dim \mathbb{U}^{(q)}$, so $\Sigma_f = \mathcal{P}_{\beta, q}$. Thus $\mathcal{G}_{\beta, 2r+1} = \mathbb{U}^{(q)} \setminus f(\mathcal{P}_{\beta, q})$ is generic and has full measure, as wanted. ∎

Lemma 2.6 is, in turn, an immediate consequence of the following fact, because $\mathcal{P}_{\beta, q}$ is an open subset of $\widetilde{\mathcal{P}}_{\beta, q}$.

**Proposition 2.7** Fix any positive integer $q$. Then, for generic $\beta \in \mathcal{C}^\infty(\mathbb{X} \times \mathbb{U}, \mathbb{R})$, the set

$$\widetilde{\mathcal{P}}_{\beta, q} := \left\{ \left((x_1, x_2), w\right) \in \mathbb{X}^{(2)} \times \mathbb{U}^{(q)} \mid \beta(x_1, w) = \beta(x_2, w) \right\}. \tag{11}$$

is either empty or it is a submanifold of $\mathbb{X}^{(2)} \times \mathbb{U}^{(q)}$ of dimension $qm + 2r - q$.

*Proof.* The proof is a routine exercise in transversality theory. We begin by recalling the Multijet Transversality Theorem (see [10] for details), in the special case of jets of order zero. (The general case for arbitrary orders, which we do not need here, would require the careful definition of jets of functions on manifolds, which in turn requires a more complicated quotient space construction.) The theorem states that, given:

- any two smooth manifolds $M$ and $N$,
- any positive integer $s$, and
- any submanifold $W$ of $M^s \times N^s$,

then, for generic $\beta \in \mathcal{C}^\infty(M, N)$, it holds that the $s$-fold 0-prolongation $j_s^0 \beta$ of $\beta$ is transversal[*] to $W$, where

$$j_s^0 \beta : \ M^{(s)} \to M^s \times N^s : \ (\xi_1, \ldots \xi_s) \mapsto (\xi_1, \ldots \xi_s, \beta(\xi_1), \ldots \beta(\xi_s)) \ .$$

Thus, since preimages of submanifolds under transversal maps are submanifolds of the same codimension, a generic $\beta \in \mathcal{C}^\infty(M, N)$ has the property that $(j_s^0 \beta)^{-1}(W)$ is either empty or it is a submanifold of $M^{(s)}$ of codimension (that is, $\dim M^{(s)} - \dim W$) equal to the codimension of $W$ (that is, $\dim M^s \times N^s - \dim W$).

A typical application of this theorem is in the context of the Whitney embedding theorem, as it implies that the set of one-to-one smooth mappings from $M$ to $N$ is generic, provided that $\dim N > 2 \dim M$. To see this, one just takes $s = 2$ and $W = \{(\xi_1, \xi_2, \zeta_1, \zeta_2) \mid \zeta_1 = \zeta_2\}$, which has codimension equal to $\dim N$. Then $\widetilde{W} = (j_s^0 \beta)^{-1}(W) \subseteq M^{(2)}$ is the set of pairs $\xi_1 \neq \xi_2$ such that $\beta(\xi_1) = \beta(\xi_2)$, and this set must be empty, since otherwise its codimension would be larger than $\dim M^{(2)} = 2 \dim M$, which is nonsense. (See e.g. [10, 14, 15, 20] for such arguments.) The application to our result is very similar, and proceeds as follows.

Pick $M = \mathbb{X} \times \mathbb{U}$, $N = \mathbb{R}$, and $s = 2q$. To define $W$, we write the coordinates in $M^s$, and in particular in the subset $M^{(s)}$, in the following form:

$$((x_1, u_1), (x_2, u_2), \ldots, (x_q, u_q), (y_1, v_1), (y_2, v_2), \ldots, (y_q, v_q)) \tag{12}$$

and let $W$ consist of those elements

$$((x_1, u_1), (x_2, u_2), \ldots, (x_q, u_q), (y_1, v_1), (y_2, v_2), \ldots, (y_q, v_q), w_1, \ldots, w_q, z_1, \ldots, z_q)$$

in $M^s \times \mathbb{R}^s$ such that:

$$x_1 = x_2 = \ldots = x_q, \ y_1 = y_2 = \ldots = y_q, \ u_1 = v_1, u_2 = v_2, \ldots, u_q = v_q, \tag{13}$$

and $w_i = z_i$ for all $i = 1, \ldots, q$. Counting equations, it is clear that $W$ is a submanifold (linear subspace in the obvious local coordinates) of codimension $\rho := 2r(q-1) + mq + q$.

For generic $\beta$, multijet transversality insures that $\mathcal{Q} = (j_s^0 \beta)^{-1}(W)$ is either empty or it is a submanifold of $M^{(s)}$ of codimension $\rho$. The set $\mathcal{Q}$ consists of all sequences (of distinct pairs) as in (12) for which all the equalities (13) hold as well as

$$\beta(x_1, u_1) = \beta(y_1, u_1), \ \beta(x_2, u_2) = \beta(y_2, u_2), \ \ldots \beta(x_q, u_q) = \beta(y_q, u_q) \ .$$

---
[*]Recall that transversality of $f : P \to Q$ to a submanifold $W$ of $Q$ denoted $f \pitchfork W$, means that $Df_x(T_x P) + T_{f(x)} W = T_{f(x)} Q$ for all $x$ such that $f(x) \in W$. All that we need here is the conclusion on preimages.

Now we introduce the function $\Phi : \mathbb{X}^{(2)} \times \mathbb{U}^{(q)} \to (\mathbb{X} \times \mathbb{U})^{(2q)}$ that maps

$$((x,y),(u_1,\ldots,u_q)) \mapsto (x,u_1),(x,u_2),\ldots,(x,u_q),(y,u_1),(y,u_2),\ldots,(y,u_q)$$

and notice that $\Phi$ establishes a diffeomorphism between $\widetilde{\mathcal{P}}_{\beta,q}$ and its image $\mathcal{Q}$. Thus, $\widetilde{\mathcal{P}}_{\beta,q}$ is either empty or is a submanifold of dimension $2q(m+r) - \rho = qm + 2r - q$, as claimed. ∎

## 3 A Special Case: Responses Linear in Inputs

A very special case of our setup is that in which the input set is an Euclidean space: $\mathbb{U} = \mathbb{R}^m$, and $\beta$ depends linearly on the inputs $u$. Then we may write:

$$\beta(x,u) \;=\; f(x) \cdot u \tag{14}$$

(dot indicates inner product in $\mathbb{U}$), where $f$ is some mapping $\mathbb{X} \to \mathbb{R}^m$. We will denote by $S \subseteq \mathbb{R}^m$ the image $f(\mathbb{X})$ of $f$. In this special case, we have that two parameters $x_1$ and $x_2$ are indistinguishable by inputs belonging to a given subset $\mathbb{U}_0 \subseteq \mathbb{R}^m$ ($x_1 \underset{\mathbb{U}_0}{\sim} x_2$) if and only if

$$f(x_1) - f(x_2) \in \mathbb{U}_0^\perp$$

($\mathbb{U}_0^\perp = \{a \in \mathbb{R}^m \mid a \cdot u = 0 \,\forall u \in \mathbb{U}_0\}$ is the orthogonal complement of $\mathbb{U}_0$), and $x_1$ and $x_2$ are indistinguishable (case $\mathbb{U}_0 = \mathbb{U}$) if and only if $f(x_1) = f(x_2)$. Therefore, a subset $\mathbb{U}_0$ is a universal distinguishing set if and only if $(S - S) \bigcap \mathbb{U}_0^\perp = \{0\}$, i.e.

$$a, b \in S, \; a - b \in \mathbb{U}_0^\perp \;\; \Rightarrow \;\; a = b\,.$$

Let us introduce the set $\sec(S)$ consisting of all unit secants of $S$:

$$\sec(S) := \left\{ \frac{a-b}{|a-b|}, \, a \neq b, \, a, b \in S \right\}$$

as well as the set $\mathrm{u}(\mathbb{U}_0^\perp)$ of unit vectors in $\mathbb{U}_0^\perp$. These are both subsets of the $(m-1)$-dimensional unit sphere $\mathbb{S}^{m-1}$. With these notations, we can say that a subset $\mathbb{U}_0$ is a universal distinguishing set if and only if $\sec(S) \bigcap \mathrm{u}(\mathbb{U}_0^\perp) = \emptyset$, or equivalently, if and only if

$$\mathrm{u}(\mathbb{U}_0^\perp) \subseteq \mathbb{S}^{m-1} \setminus \sec(S)\,. \tag{15}$$

Any basis of $\mathbb{U}$ provides a universal distinguishing set of cardinality $m$ (since then $\mathrm{u}(\mathbb{U}_0^\perp) = \emptyset$). On the other hand:

**Proposition 3.1** There is a universal distinguishing set of cardinality $m - 1$ if and only if $\sec(S) \neq \mathbb{S}^{m-1}$.

*Proof.* If there is such a $\mathbb{U}_0$ with less than $m$ elements, then $\mathbb{U}_0^\perp \neq \{0\}$, so also $\mathrm{u}(\mathbb{U}_0^\perp) \neq \emptyset$, and then (15) gives that $\sec(S) \neq \mathbb{S}^{m-1}$. Conversely, if there is any $u \in \mathbb{S}^{m-1} \setminus \sec(S)$, we may let $\mathbb{U}_0$ be any basis of $\{u\}^\perp$, so that $\mathrm{u}(\mathbb{U}_0^\perp) = \{\pm u\}$. Since also $-u \in \mathbb{S}^{m-1} \setminus \sec(S)$, we have that $\mathrm{u}(\mathbb{U}_0^\perp) \subseteq \mathbb{S}^{m-1} \setminus \sec(S)$, and thus $\mathbb{U}_0$ is a universal distinguishing set of cardinality $m - 1$. ∎

### 3.1 Examples

We will provide examples of two subsets $\mathcal{S} \subseteq \mathbb{R}^2$ and $\mathcal{R} \subseteq \mathbb{R}^3$, both images of analytic maps defined on $\mathbb{R}$, with the properties that $\mathcal{S} - \mathcal{S} = \mathbb{R}^2$ and $\sec(\mathcal{R}) = \mathbb{S}^2$. These examples will allow us to show that the $2r+1$ bound is best possible. We need this first:

**Lemma 3.2** Pick any real $a > 0$ and any nonnegative integer $k$, and let

$$f(x) = (x+a)\sin x \tag{16}$$

for $x \in \mathbb{R}$. Then, there exists an

$$M_k \in ((2k+1/2)\pi, (2k+1)\pi)$$

and a continuous map

$$\alpha : [2k\pi, M_k] \to [M_k, (2k+1)\pi]$$

such that

$$\alpha(2k\pi) = (2k+1)\pi, \ \alpha(M_k) = M_k, \ \text{and} \ f(\alpha(x)) = f(x) \, \forall \, x \in [2k\pi, M_k].$$

*Proof.* Consider the restriction of the function $f$ to the interval $[2k\pi, (2k+1)\pi]$, and observe that its derivative $f'(x) = \sin x + (x+a)\cos x$ is positive for $2k\pi \leq x \leq (2k+1/2)\pi$, has $f'((2k+1)\pi) = ((2k+1)\pi + a) \cdot (-1) < 0$, and, for $(2k+1/2)\pi < x < (2k+1)\pi$, $f'(x) = 0$ is equivalent to

$$\tan x = -x - a,$$

which happens at a unique $x = M_k \in ((2k+1/2)\pi, (2k+1)\pi)$ (clear from the graph of $\tan x$ and from the fact that the graph of $-x-a$ is in the fourth quadrant). Therefore, on the interval $[2k\pi, (2k+1)\pi]$, $f$ is strictly increasing on $[2k\pi, M_k]$ and strictly decreasing on $[M_k, (2k+1)\pi]$. Let $f_1$ and $f_2$ be the restrictions of $f$ to $[2k\pi, M_k]$ and $[M_k, (2k+1)\pi]$ respectively, and let

$$g := f_2^{-1} : [0, f(M_k)] \to [M_k, (2k+1)\pi]$$

(so $g$ is a strictly decreasing continuous function). Finally, let $\alpha := g \circ f_1$. Thus $\alpha$ is a continuous function and it satisfies that $\alpha(2k\pi) = (2k+1)\pi$ and $\alpha(M_k) = M_k$ by construction. Finally,

$$f(\alpha(x)) = f_2(\alpha(x)) = f_2(f_2^{-1}(f_1(x))) = f_1(x) = f(x)$$

for all $x \in [2k\pi, M_k]$, as desired. ■

**Lemma 3.3** Consider the spiral $\mathcal{S} = \{\zeta \in \mathbb{C} \mid \zeta = re^{ir}, r \geq 0\}$. Then, for each complex $z \in \mathbb{C}$ there exist two elements $\zeta_1, \zeta_2 \in \mathcal{S}$ such that $\zeta_1 - \zeta_2 = z$. That is, as a subset of $\mathbb{R}^2$, we have $\mathcal{S} - \mathcal{S} = \mathbb{R}^2$.

*Proof.* The idea of the proof is very simple: we first find some "chord" between two points $a$ and $b$ such that the difference $b - a$ is a multiple of the desired $z$, and its modulus is larger than that of $z$; then, we displace it in an orthogonal direction, until the resulting chord has the right length. Analytically, we proceed as follows.

Let $z = re^{i\varphi}$, with $r \geq 0$ and $\varphi \geq 0$. Without loss of generality, we assume that $r > 0$ (if $z = 0$, we just take $z_1 = z_2$ to be any element of $\mathcal{S}$), and $\varphi > 0$, and pick any positive integer $k$ such that $\varphi + 2k\pi > \frac{r}{2}$. Our goal is, thus, to show that there are reals $s, t$ such that

$$te^{i(t-\varphi)} - se^{i(s-\varphi)} = r$$

which is equivalent to asking that

$$t \sin(t - \varphi) = s \sin(s - \varphi) \tag{17}$$

and

$$t \cos(t - \varphi) - s \cos(s - \varphi) = r. \tag{18}$$

By Lemma 3.2, applied with $a = \varphi$ and the chosen $k$, there exists $M_k \in ((2k+1/2)\pi, (2k+1)\pi)$ and a continuous $\alpha : [2k\pi, M_k] \to [M_k, (2k+1)\pi]$ such that $\alpha(2k\pi) = (2k+1)\pi$, $\alpha(M_k) = M_k$, and

$$(x + \varphi) \sin x = (\alpha(x) + \varphi) \sin \alpha(x) \tag{19}$$

for all $x \in [2k\pi, M_k]$. Let

$$\gamma(t) := t \cos(t - \varphi) - [\varphi + \alpha(t - \varphi)] \cos(\alpha(t - \varphi)), \ t \in [\varphi + 2k\pi, \varphi + M_k].$$

Then,

$$\gamma(\varphi + 2k\pi) = (\varphi + 2k\pi) \cos(2k\pi) - (\varphi + 2k\pi + \pi) \cos(2k\pi + \pi) = 2(\varphi + 2k\pi) + \pi > r$$

and

$$\gamma(\varphi + M_k) = (\varphi + M_k) \cos M_k - (\varphi + M_k) \cos M_k = 0.$$

So, since $\gamma$ is continuous, there is some $t_0 \in [\varphi + 2k\pi, \varphi + M_k]$ such that $\gamma(t_0) = r$, which means that

$$t_0 \cos(t_0 - \varphi) - s_0 \cos(s_0 - \varphi) = r$$

with $s_0 := \varphi + \alpha(t_0 - \varphi)$. Moreover, evaluating (19) at $x = t_0 - \phi$ gives that

$$t_0 \sin(t_0 - \varphi) = s_0 \sin(s_0 - \varphi),$$

and the proof is now complete. ∎

**Lemma 3.4** Consider the following subset of $\mathbb{R}^3$:

$$\mathcal{R} = \{(\zeta, \xi) \in \mathbb{C} \times \mathbb{R} \mid \zeta = e^{ir}, \ \xi = e^r \sin 2r, \ r \geq 0\}.$$

Then, for each $w \in \mathbb{C} \times \mathbb{R}$ there exist two elements $\omega_1, \omega_2 \in \mathcal{R}$ and a real number $\chi > 0$ such that $\omega_1 - \omega_2 = \chi w$. That is, as a subset of $\mathbb{R}^3$, $\sec(\mathcal{R}) = \mathbb{S}^2$.

*Proof.* We will use the real functions

$$g(x) := e^x \sin 2x \ \text{ and } \ f(x) := g(x - \pi/2).$$

Representing complex numbers $\zeta$ in the form $re^{i\varphi}$, with $r$ not necessarily positive but the argument $\varphi$ restricted to $[0, \pi)$, we may rephrase the claim of the lemma as follows: for each $r \in \mathbb{R}$, $0 \leq \varphi < \pi$, and $p \in \mathbb{R}$, there exist $\chi > 0$ and $t, s \geq 0$ such that

$$\begin{pmatrix} e^{it} - e^{st} \\ -e^{\pi/2}(f(t) - f(s)) \end{pmatrix} = \chi \begin{pmatrix} re^{i\varphi} \\ p \end{pmatrix} \tag{20}$$

and we may assume that $w = (r, p) \neq 0$, since when $w = 0$ we may pick any $\chi$ and $\omega_1 = \omega_2 =$ any element of $\mathcal{R}$. Equality of the first (complex) components in (20) amounts to asking that the following two equations hold:

$$\sin(t - \varphi) = \sin(s - \varphi), \quad \cos(t - \varphi) - \cos(s - \varphi) = \chi r, \tag{21}$$

and for this, in turn, it is sufficient to find $t, s \in \mathbb{R}$ such that:

$$t + s = \pi + 2\varphi, \quad \cos(\pi + \varphi - s) - \cos(s - \varphi) = \chi r.$$

The last of these equations simplifies to $-\cos(s - \varphi) = \chi r/2$. In summary, for any given $r, \varphi, p$ we must find two real numbers $\alpha$ and $s$ such that (absorbing $-e^{-\pi/2}$ into $p$):

$$f(\pi + 2\varphi - s) - f(s) = p, \quad -\cos(s - \varphi) = \chi r/2.$$

Or, letting $\theta := s - \varphi - \pi/2$, and in terms of $g(x) = f(x + \pi/2)$, the problem is to solve the following two simultaneous equations for $\chi, \theta$:

$$g(\varphi + \theta) - g(\varphi - \theta) = \chi p, \quad \sin \theta = \chi r/2. \tag{22}$$

In the special case that $r = 0$, we have $p \neq 0$ (recall $(r, p) \neq 0$), so we may pick $\theta := \pi$ and $\chi := [g(\varphi + \pi) - g(\varphi - \pi)]/p$. Thus, we assume $r \neq 0$ from now on.

We will show that there is some $\theta$ which is not of the form $k\pi$ for any integer $k$, such that

$$\frac{g(\varphi + \theta) - g(\varphi - \theta)}{\sin \theta} = \frac{2p}{r}. \tag{23}$$

Once such a $\theta$ is found, we may simply let $\chi := (2/r) \sin \theta$, from which it follows that $\sin \theta = \chi r/2$ and $g(\phi + \theta) - g(\varphi - \theta) = (2/r) p \sin \theta = \chi p$, so (22) holds as desired.

Letting $q := (2p/r) e^{-\varphi}$, we restate our goal as that of solving

$$\frac{e^\theta \sin(2(\varphi + \theta)) - e^{-\theta} \sin(2(\varphi - \theta))}{\sin \theta} = q \tag{24}$$

for $q$. The idea is to take $\theta \gg 0$, so that the first term in the numerator dominates. We consider three cases: (i) $0 < \varphi < \pi/2$, (ii) $\varphi = 0$, and (iii) $\pi/2 \leq \varphi < \pi$. The last case follows from the first two, since given any $\varphi$ in the interval $[\pi/2, \pi)$, we may solve the version of (24) stated for $\hat{\varphi} := \varphi - \pi/2$ instead of $\varphi$ and $-q$ instead of $q$, and the same $\theta$ then solves (24) (since $\sin(2(\hat{\varphi} + \theta)) = -\sin(2(\varphi + \theta))$ and $\sin(2(\hat{\varphi} - \theta)) = -\sin(2(\varphi - \theta))$).

*Case (i):* We introduce the functions $A(\theta) := e^\theta \sin(2(\varphi + \theta))$ and $B(\theta) := e^{-\theta} \sin(2(\varphi - \theta))$, and note that $|B(\theta)| < 1$ for $\theta > 0$. Pick

$$\alpha := \frac{1}{4} \min\left\{\varphi, \frac{\pi}{2} - \varphi\right\} > 0$$

and observe that

$$4\alpha \leq \varphi \leq \pi/2 - 4\alpha \tag{25}$$

(from which it also follows that $\alpha \leq \pi/16$). Let

$$\sigma_1 := e^{2\alpha - \pi/2} \sin(2\varphi + 4\alpha - \pi), \quad \sigma_2 := e^{-2\alpha} \sin(2\varphi - 4\alpha).$$

Observe that $\varphi \leq \pi/2 - 4\alpha$ implies that

$$-\pi < 2\varphi + 4\alpha - \pi < -4\alpha < 0$$

so $\sigma_1 < 0$, and that $4\alpha \leq \varphi$ implies

$$0 < 4\alpha < 2\varphi - 4\alpha < 2\varphi < \pi$$

so $\sigma_2 > 0$. Now pick an odd integer $k$ large enough so that

$$\frac{e^{k\pi}\sigma_1 + 1}{\cos 2\alpha} < q \quad \text{and} \quad \frac{e^{k\pi}\sigma_2 - 1}{\sin 2\alpha} > q$$

and introduce $\theta_1 := k\pi - \pi/2 + 2\alpha$, $\theta_2 := k\pi - 2\alpha$, and the interval $I := [\theta_1, \theta_2]$. On this interval, $\sin\theta$ is decreasing and positive; in fact it satisfies $\sin 2\alpha = \sin\theta_2 \leq \sin\theta \leq \sin\theta_1 = \cos 2\alpha$. In particular, the function

$$C(\theta) := \frac{A(\theta) + B(\theta)}{\sin\theta}$$

is well-defined and continuous on $I$. Moreover, since $A(\theta_1) = e^{k\pi}\sigma_1$, $A(\theta_2) = e^{k\pi}\sigma_2$, and $|B(\theta)| < 1$ for all $\theta$, and because of the choice of $k$,

$$C(\theta_1) < q \quad \text{and} \quad C(\theta_2) > q$$

so we conclude that, indeed, we can solve $C(\theta) = q$.

*Case (ii):* If $\varphi = 0$, then (24) reduces to $\frac{e^\theta \sin 2\theta - e^{-\theta} \sin 2\theta}{\sin\theta} = q$ or equivalently

$$2e^\theta \cos\theta - 2e^{-\theta}\cos\theta = q.$$

We pick a positive integer $k$ large enough so that

$$-\sqrt{2}e^{2k\pi + 3\pi/4} + 2 < q \quad \text{and} \quad \sqrt{2}e^{2k\pi + \pi/4} - 2 > q.$$

Now consider the function $C(\theta) = 2e^\theta \cos\theta - 2e^{-\theta}\cos\theta$ on the interval $I = [\theta_1, \theta_2] = [2k\pi + \pi/4, 2k\pi + 3\pi/4]$. We have that $C(\theta_1) > q$ and $C(\theta_2) < q$, so we can again solve $C(\theta) = q$. ∎

**Corollary 3.5** For any fixed positive integer $r$, consider the subset $\mathcal{R}_r = \mathcal{S}^{r-1} \times \mathcal{R}$ of $\mathbb{R}^{2r+1}$. Then $\sec(\mathcal{R}_r) = \mathbb{S}^{2r}$.

*Proof.* Take any $\theta \in \mathbb{S}^{2r}$, and write it in the form $(z_1, \ldots, z_{r-1}, w)$ with $z_i \in \mathbb{C}$ and $w \in \mathbb{C} \times \mathbb{R}$. Using Lemma 3.4, we pick $\omega_1, \omega_2 \in \mathcal{R}$ and $\chi > 0$ such that $\omega_1 - \omega_2 = \chi w$. Next, using Lemma 3.3, we find $\zeta_{ij} \in \mathcal{S}$, $i = 1, \ldots, r-1$, $j = 1, 2$, such that $\zeta_{i1} - \zeta_{i2} = \chi z_i$ for each $i = 1, \ldots, r-1$. So $a_j := (\zeta_{1j}, \ldots, \zeta_{r-1,j}, \omega_j) \in \mathcal{R}_r$ for $j = 1, 2$ satisfy $a_1 - a_2 = \chi\theta$. Since $\theta$ has unit norm, it follows that $\chi = |a_1 - a_2|$, so $\theta = \frac{a_1 - a_2}{|a_1 - a_2|} \in \sec(\mathcal{R}_r)$. ∎

## 3.2 The Bound $2r+1$ is Best Possible

We now present an example which shows that the number $2r+1$ in Theorem 2 cannot be lowered. For this, we must exhibit, for each positive integer $r$, an analytic response $\beta$ with the property that $\mathcal{G}_{\beta,2r}$ either is not generic or has less than full measure in $\mathbb{U}^{2r}$. In fact, we will show far more: we will show that $\mathcal{G}_{\beta,2r}$ is empty.

The example is as follows. Given any fixed $r$, we consider the mapping

$$g \,:\, [0,\infty)^{r-1} \times (0,\infty) \to \mathbb{R}^{2r+1} \,:\, (t_1,\ldots,t_r) \mapsto \left(t_1 e^{it_1},\ldots,t_{r-1}e^{it_{r-1}}, e^{it_r}, e^{t_r}\sin 2t_r\right)$$

whose image is $\mathcal{R}_r$ (note that $(1,0,0)$ can be obtained as $(e^{it_r}, e^{t_r}\sin 2t_r)$ for $t_r = 2\pi$, so $t_r = 0$ is not required), and let $f : \mathbb{X} = \mathbb{R}^{r-1} \times (0,\infty) \to \mathbb{U} = \mathbb{R}^{2r+1}$ be given by $f(t_1,\ldots,t_r) = g(t_1^2,\ldots,t_{r-1}^2,t_r)$, so also $f(\mathbb{X}) = \mathcal{R}_r$. We let $\beta(x,u) = f(x) \cdot u$. By Proposition 3.1, if there is a universal distinguishing set of cardinality $2r$ then $\sec(\mathcal{R}_r) \neq \mathbb{S}^{2r}$ (note that $m-1 = 2r$). This contradicts Corollary 3.5, so no such set can exist. ∎

We can modify this example so that the input set $\mathbb{U}$ is scalar, as follows. Let us consider the following response, with $\widetilde{\mathbb{U}} = \mathbb{R}$:

$$\widetilde{\beta}(x,u) \,:=\, \beta(x,\psi(u)) \,=\, f(x) \cdot \psi(u)$$

where $\psi : \mathbb{R} \to \mathbb{R}^{2r+1} : u \mapsto (1,u,u^2,\ldots,u^{2r})$, leaving $f$ and $\mathbb{X}$ unchanged. We claim that there is no universal distinguishing set of cardinality $2r$. Indeed, suppose that $\widetilde{\mathbb{U}}_0$ would be a $2r$-element universal distinguishing set. Consider the set $\mathbb{U}_0 := \psi(\widetilde{\mathbb{U}}_0)$. As this set has $2r$ elements, it cannot be a universal distinguishing set for $\beta$. Thus, there exist two parameters $x_1$ and $x_2$ which are distinguishable for $\beta$, that is $f(x_1) \neq f(x_2)$, but such that $\beta(x_1,v) = \beta(x_2,v)$ for all $v \in \mathbb{U}_0$, which implies $\widetilde{\beta}(x_1,u) = \widetilde{\beta}(x_2,u)$ for all $u \in \widetilde{\mathbb{U}}_0$. If we show that $\widetilde{\beta}(x_1,u) \neq \widetilde{\beta}(x_2,u)$ for some $u \in \widetilde{\mathbb{U}}$ then we will have a contradiction with $\widetilde{\mathbb{U}}_0$ being a universal distinguishing set. To see this, simply notice that $\psi(\mathbb{R})$ linearly spans $\mathbb{R}^{2r+1}$: if the vector $a \in \mathbb{R}^{2r+1}$ is nonzero, then $a \cdot \psi(u) \neq 0$ for some $u \in \mathbb{R}$ ($\sum a_i u^i \equiv 0 \Rightarrow a_i = 0\, \forall i$), and apply with $a = f(x_1) - f(x_2)$. ∎

## 4 Application to Systems

Now we apply the results about abstract responses to the special case of identifying parameters in systems, proving Theorem 1. We take an analytically parametrized system $\Sigma$, with $r = \dim \mathbb{X}$. When the number of measurements $p = 1$, the results follows from Theorem 2 applied to $\beta = \beta_\Sigma$.

For the general case, we consider the scalar responses $\beta_\Sigma^i$, $i = 1,\ldots,p$ which are obtained as coordinate projections of $\beta_\Sigma$. We claim that, for each fixed $q$, and with the obvious notations, $\bigcap_i \mathcal{G}_q^i \subseteq \mathcal{G}_q$. Indeed, take any $w \in \bigcap_i \mathcal{G}_q^i$, and any $x_1, x_2$ such that $x_1 \not\sim x_2$. Then there must be some $i \in \{1,\ldots,q\}$ such that $x_1 \not\sim x_2$ for the response $\beta_\Sigma^i$, so, since $w \in \mathcal{G}_q^i$, it follows that $\beta_\Sigma^i(x_1,w) \neq \beta_\Sigma^i(x_2,w)$, and therefore also $\beta_\Sigma(x_1,w) \neq \beta_\Sigma(x_2,w)$. This proves that $w \in \mathcal{G}_q$. Since the intersection of a finite (or even countable) number of generic and full measure sets is again generic and of full measure, $\mathcal{G}_q$ must have this property. This completes the proof of Theorem 1 for arbitrary $p$. ∎

### 4.1 A System for which $2r+1$ Experiments are Needed

We can express the responses $\beta$ or $\widetilde{\beta}$ from Section 3.2 as the response $\beta_\Sigma$ for an analytically parametrized system, and in this way know that the number $2r+1$ in Theorem 1 cannot be lowered to $2r$, and in fact, there are analytic systems with $r$ parameters for which there is not even a single universal distinguishing set of cardinality $2r$. The simplest $\Sigma$ would be obtained by using any $f$, and just defining $h(z, u, x) = \beta(x, u)$. It is far more interesting, however, to give an example where only polynomials appear in the system description and $h$ does not depend directly on $u$ and $x$. We do this explicitly for the case $r = 1$; the case of arbitrary $r$ is entirely analogous.

The system $\Sigma$ that we construct has state space $M = \mathbb{R}^9$, input-value space $U = \mathbb{R}$, parameter space $\mathbb{X} = (0, \infty)$, experiment space $\Lambda = \mathbb{R}$, and $p = 1$, and is given by:

$$\chi \;:\; (0,\infty) \to \mathbb{R}^9 \;:\; a \mapsto (a, 1, 0, 0, 1, 0, 1, 0, 1),$$

$$h \;:\; \mathbb{R}^9 \to \mathbb{R} \;:\; (z_1, z_2, \ldots, z_9) \mapsto z_2 z_5 + z_3 z_6 + z_4 z_8 z_9$$

(independent of $u$ and $x$),

$$\mu \;:\; \mathbb{R} \times \mathbb{R} \to \mathbb{R} \;:\; (\lambda, t) \mapsto \lambda$$

(i.e., inputs are constant scalars), and

$$f(z, u, x) \;:=\; (0, 0, u z_2, 2 u z_3, -z_1 z_6, z_1 z_5, -2 z_1 z_8, 2 z_1 z_7, z_1 z_9).$$

The solution $z(t)$ with initial condition $z(0) = \chi(a)$ and input $u(t) \equiv \lambda$ is the following vector:

$$(a, 1, \lambda t, \lambda^2 t^2, \cos at, \sin at, \cos 2at, \sin 2at, e^{at})$$

and therefore

$$\beta_\Sigma(a, \lambda) = h(z(1)) = \varphi(a) \cdot \psi(\lambda)$$

where $\varphi(a) = (\cos a, \sin a, e^a \sin 2a)$ and $\psi(\lambda) = (1, \lambda, \lambda^2)$, which is $\widetilde{\beta}(a, \lambda)$. ∎

A small modification of this example has $h$ linear: just add an additional variable $z_{10}$ with initial condition $z_{10}(0) = 0$ and satisfying $\dot{z}_{10} = (z_2 z_5 + z_3 z_6 + z_4 z_8 z_9)^{\cdot}$ (written, using $\dot{z}_3 = u z_2$, etc, in terms of the $z_i$ and $u$), and now use $h(z) = z_{10}$.

## 5 Distinguishability in the Operon Example

We now show, for the operon system (5) with external input, that every two distinct parameters are distinguishable. We work out such an example in order to emphasize that the problem of determining identifiability is nontrivial, which makes it more interesting that Theorem 1 applies without this knowledge. The experiments $(\lambda, T)$ consist of using constant inputs $u(t) \equiv \lambda$ for varying intervals $[0, T]$ and measuring $M(T)$. Thus, we wish to prove: if for every nonnegative $\lambda$, the solution $(M(t), E(t))$ of $\dot{M} = E^m/(1 + E^m) - aM$, $\dot{E} = M - bE - \lambda E$, with initial condition $(M_0, E_0)$, and the solution $(M^\dagger(t), E^\dagger(t))$ of $\dot{M}^\dagger = (E^\dagger)^{m^\dagger}/(1 + (E^\dagger)^{m^\dagger}) - a^\dagger M^\dagger$, $\dot{E}^\dagger = M^\dagger - b^\dagger E^\dagger - \lambda E^\dagger$, with initial condition $(M^\dagger{}_0, E^\dagger{}_0)$ are such that $M(t) \equiv M^\dagger(t)$, then necessarily $M_0 = M^\dagger{}_0$, $E_0 = E^\dagger{}_0$, $m = m^\dagger$, $a = a^\dagger$, and $b = b^\dagger$.

Assume given two parameters $(M_0, E_0, m, a, b)$ and $(M^\dagger{}_0, E^\dagger{}_0, m^\dagger, a^\dagger, b^\dagger)$ with this property (recall that the entries are all positive, by assumption). Since $M(t) \equiv M^\dagger(t)$, of course $M_0 = M^\dagger{}_0$, and we write $\xi$ for their common value. Now fix an arbitrary $\lambda$ and look at $M(1)$ and $M^\dagger(1)$. We have that $M(1) = e^{-a}\xi + \int_0^1 e^{-a(1-t)} \alpha(t)\, dt$, where we define $\alpha(t) = \frac{E(t)^m}{1+E(t)^m}$ and note that $\alpha(t) \leq 1$ for all $t$. It follows that $M(1) \leq 1 + \xi$ is bounded independently of the value of $\lambda$. On the other hand, $E(t) = e^{-(b+\lambda)t}E_0 + \int_0^t e^{-(b+\lambda)(t-s)}M(s)\,ds$. Since $e^{-(b+\lambda)(t-s)} \to 0$ as $\lambda \to \infty$ for each $s < t$, and $M$ is bounded, we have by dominated convergence that $E(t) \to 0$ as $\lambda \to \infty$, for each fixed $t$. Thus also $\alpha(t) \to 0$ as $\lambda \to \infty$, for each fixed $t$. Now applying this to the above formula for $M(1)$, and again by dominated convergence, we have that $M(1) \to e^{-a}\xi$ as $\lambda \to \infty$. Since $M(1) = M^\dagger(1)$ for any given $\lambda$, and since by an analogous argument we also have that $M^\dagger(1) \to e^{-a^\dagger}\xi$ as $\lambda \to \infty$, we conclude that $e^{-a}\xi = e^{-a^\dagger}\xi$, and therefore that $a = a^\dagger$.

From the original differential equation $\dot{M} = \frac{E^m}{1+E^m} - aM$, we know that $\alpha = \dot{M} + aM$, and we also have the same formula ($a$ is the same, and $M$ too) for the second set of parameters, which gives us that $\alpha(t) = \alpha^\dagger(t)$, and therefore (since $p \to \frac{p}{1+p}$ is one to one) that $E^m(t) = (E^\dagger)^{m^\dagger}(t)$ for all $t$, no matter what input $\lambda$ is used. For any given $\lambda$, we introduce the function

$$w(t) = \frac{1}{\lambda}\frac{1}{E^m(t)}\frac{d}{dt}E^m(t)$$

and similarly for the second set of parameters. As $E^m \equiv (E^\dagger)^{m^\dagger}$, also $w \equiv w^\dagger$. Calculating, we have that $w = \frac{m}{\lambda E}(M - bE) - m$, and similarly for $w^\dagger$. Thus we obtain, evaluating at $t = 0$:

$$\frac{m}{\lambda E_0}(\xi - bE_0) - m \equiv \frac{m^\dagger}{\lambda E^\dagger_0}(\xi - b^\dagger E^\dagger_0) - m^\dagger$$

and taking now the limit as $\lambda \to \infty$ we conclude that $m = m^\dagger$. Thus, since now we know that $E^m \equiv (E^\dagger)^m$, we can conclude that $E \equiv E^\dagger$ and in particular that $E_0 = E^\dagger_0$ and that $dE/dt \equiv dE^\dagger/dt$ (for any given value of the input $\lambda$). Finally, for $\lambda = 0$ we have from $\dot{E} = M - bE$ that $b = b^\dagger$. ∎

# 6 Comments and Relations to Other Work

We close with some general comments.

## 6.1 Universal Distinguishing Sets

The concepts of distinguishability and distinguishing sets are common in several fields. In control theory (see e.g. [25], Chapter 6), one studies the possibility of separating internal states (corresponding to parameters in the current context) on the basis of input/output experiments. The papers [23, 24] deal with applications of distinguishability to the study of local minima of least-squares error functions, and set-shattering in the sense of Vapnik and Chervonenkis, for artificial neural networks. More combinatorial, but essentially the same, notions, have appeared in computational learning theory (a *teaching set* is one which allows a "teacher" to uniquely specify the particular function being "taught" among all other functions of interest, see e.g. [9]) and in the theory of experiments in automata and sequential machine theory (cf. [5]).

## 6.2 Observability

The observability problem, that is, the reconstruction of all internal states of the system, is included in the problem discussed here, in the following sense. Suppose that parameters include all initial states, that is, $\mathbb{X} = M \times \mathbb{X}_0$ and the initial state $\chi$ is a projection onto the components in $M$ (as in the examples in Section 1.1). Then distinguishability of parameters implies distinguishability of initial states, and, since the flow of a differential equation induces a group of diffeomorphisms (so the map $z(0) \mapsto z(T)$ is one-to-one, for each $T$ and each given input), also the distinguishability of states at any future time.

## 6.3 Restarting

Notice an important feature of the setup. Since the objective is to find parameters, and these are constant, it is implicitly assumed that one may "restart" different experiments at the same initial state. In practice, this may or may not be a valid assumption. In fact, much work in control theory deals with identification problems for which one need not restart the system: this is the subject of the area of *universal inputs* for observability, cf. [12, 21, 26, 30]. On the other hand, in molecular biology multiple experiments, assuming identical initial conditions, are usually performed by careful assay controls, or by dealing with synchronized daughter cells. Indeed, because of the noisiness inherent in biological applications, data for a "single experiment" may actually represent an average of different runs under the same (approximate) conditions. In addition, many measurements in cell biology are destructive, and thus is impossible to take measurements at different times from the same cell, so the theory of universal inputs does not apply under such circumstances.

## 6.4 Genericity

The material in Section 2.3 on genericity is motivated by, and shares many of the techniques with, the theory of manifold embeddings (see also Section 6.6 below, as well as the remark in the proof about one to one maps). Closely related is also the work of Takens [28], which shows that generically, a smooth dynamical system on an $r$-dimensional manifold can be embedded in $\mathbb{R}^{2r+1}$, as well as the control-theory work of Aeyels on generic observability, which shows in [2] that for generic vector fields and observation maps on an $r$-dimensional manifold, $2r+1$ observations at randomly chosen times are enough for observability, and in [1] that this bound is best possible. Aeyels proofs, in particular, are based on transversality arguments of the general type that we use.

## 6.5 The Examples

In Lemma 3.4, if instead of $\mathcal{R}$ we would have considered the set consisting of those $(\zeta, \xi) \in \mathbb{C} \times \mathbb{R}$ with $\zeta = e^{ir}$ and $\xi = \sin 2r$, then $\sec(\mathcal{R})$ would be a proper subset of $\mathbb{S}^2$. Indeed, this set was studied in [4], where it was shown to have nonzero measure and a complement also of nonzero measure.

## 6.6 Whitney's Embedding Theorem

The material in Section 3 is closely related to the proof of the "easy version" of Whitney's Embedding Theorem, cf. [10, 14]. We briefly review this connection here.

We suppose that $S$ is a compact $r$-dimensional embedded submanifold of $\mathbb{U} = \mathbb{R}^m$. We assume that $m \geq 2r+1$ (otherwise what follows is not interesting.) The dimension $r$ may well be smaller than the dimension of $\mathbb{X}$. Of course, there is no reason whatsoever for $S$ to be a submanifold of $\mathbb{R}^m$, as the mapping $f$ may well have singularities. Thus, we are imposing yet another condition besides linearity on $u$. On the other hand, analyticity of $f$, i.e. analyticity of $\beta$ on $x$, is not needed in what follows.

The facts that there exist universal distinguishing sets $\mathbb{U}_0$ of cardinality $2r+1$ and, moreover, that a $\mathbb{U}_0$ has this property are, in the special case being considered here, an immediate consequence of the proof of the "easy version" of Whitney's Embedding Theorem. (The classical embedding results date back to Menger's 1926 work (cf. [16]) for continuous functions and maps from topological spaces into Euclidean spaces, and the smooth version dealing with embeddings of differentiable manifolds of dimension $r$ in $\mathbb{R}^{2r+1}$ due to Whitney in [31]. A "harder" version of Whitney's theorem (cf. [32]) shows that one may embed such manifolds in $\mathbb{R}^{2r}$ as well, and locally embed (immerse) in $\mathbb{R}^{2r-1}$ when $r > 1$, see [14].)

Let us briefly sketch how this conclusion is obtained. We consider first the special case $m = 2r+2$. The universal distinguishing sets consisting of $2r+1$ linearly independent vectors are in a one-to-one correspondence (up to a choice of basis) with the possible $2r+1$-dimensional subspaces $V$ of $\mathbb{R}^m$ for which $\sec(S) \bigcap u(V^\perp) = \emptyset$, or equivalently with the unit vectors $u \in \mathbb{R}^m$ which do not belong to $\sec(S)$. (Note that $u \in \sec(S)$ if and only if $-u \in \sec(S)$, so there is no need to work with projective space in dimension $m-1$, and we may simply deal with unit vectors.) Thus, one needs to show that $\mathbb{S}^{m-1} \setminus \sec(S)$ is generic and of full measure. Now, $\sec(S)$ is the image of the differentiable mapping

$$\widetilde{S} = \{(a,b) \in S^2 \mid a \neq b\} \to \mathbb{S}^{m-1} : (a,b) \mapsto \frac{a-b}{|a-b|}$$

and $\widetilde{S}$ has dimension $2r < m - 1 =$ dimension of $\mathbb{S}^{m-1}$. Thus, the Morse-Sard Theorem says that $\sec(S)$ has measure zero and is included in a countable union of closed nowhere dense sets, as wanted. The general case ($m \geq 2r + 3$) can be obtained inductively, by iteratively reducing to a smaller dimensional embedding space by means of projections along vectors $u$ picked as in the previous discussion, with a small modification: the choice of $u$ has to be made with some care, requiring in addition that all tangents to $S$ miss $u$; when doing so, the projection of $S$ has a manifold structure and the argument can indeed be repeated. See [14] for details, and also [4] for an expository discussion of these ideas in the context of numerical algorithms which optimize the projections; the generalization of the material in this last reference, to cover special classes of nonlinear parametrizations, would be of great interest.

## 6.7 The Techniques

As we mentioned, the main result is based on the facts regarding analytic functions which we developed in our previous paper [23]. This is in contrast to work based on Whitney embeddings and transversality arguments. Quite related, on the other hand, is the recent (and independent of [23]) work [7, 18], which deals with the distinguishability of fluid flows on the basis of a finite set of exact experiments: the authors provide a bound of the form "$16r + 1$" measurements (the number arises from the need to obtain appropriate parametrizations), where $r$ is the dimension of an attractor for the systems being studied, and they also employ analytic-function techniques in a manner very similar to ours.

## 6.8 Least Squares

We have not discussed the actual numerical computation of parameters on the basis of experiments, which is of course a most important direction of study. What we can say is the following: if all distinct parameters are distinguishable, then, for a generic set of $2r+1$ experiments, global minimization of the least-squares fit error function will result in a unique global minimum. But nothing is said about local minima nor, certainly, the effect of noise. (To study such effects, one will have to combine techniques as here with classical statistical tools, such as the Cramer-Rao inequality for the Fisher information matrix, which lower-bounds the covariance of any unbiased estimator. However, our global results are in any case of a rather different nature than these classical statistical techniques, which are closer to the nonsingular case treated in Section 2.2.) Note that even if not all parameters are distinguishable, the results in this paper might still be useful, see [23] for related work.

## 6.9 Vector Outputs

The statement of the main theorems notwithstanding, the results are really about scalar measurements, in the sense that the number of simultaneous measurements $p$ does not enter into the estimates. This is unavoidable: for instance, if all coordinates of $h$ happen to be the same, no additional information can be gained. It would be of interest to come up with a natural (and verifiable) condition of independence which, when incorporated into the system description, would allow one to introduce a factor $1/p$ into the estimates. It is fairly obvious how to do such a thing with abstract responses and if there are enough input dimensions ($m \geq p$): provided that independence implies that the codimensions of the sets $\mathcal{U}(x_1, x_2)$ is $p$ instead of 1, the critical inequality $(2r+1)(m-1) + 2r < (2r+1)m$ becomes $(\frac{2r}{p}+1)(m-p) + 2r < (\frac{2r}{p}+1)m$, so $\frac{2r}{p}+1$ randomly chosen experiments suffice. But in the case of systems, and even for abstract responses with low-dimensional $\mathbb{U}$, how to state a good result is less clear.

## 6.10 Structure

The problem of structure determination, that is to say, finding the *form* of equations, can sometimes be reduced to the problem studied here. Specifically, it usually happens in applications that one merely wishes to know if a particular term appears or not in the description of a differential equation. As an illustration, take the following situation in molecular biology: it is not known if a variable, let us say $z_1$, affects or not the evolution of another variable, let us say $z_2$, but it is known that, *if* there is any effect at all, then this influence takes the form of an inhibitory feedback term $c\frac{1}{1+z_2}z_1$ appearing in the equation for $\dot{z}_1$. One reduces to the identification problem by thinking of "$c$" as a parameter; $c = 0$ corresponds to no effect. Given that the "hypothesis testing" problem "determine if $c = 0$ or $c \neq 0$" is less demanding than the problem of actually finding the value of $c$, it is not surprising that less than $2r+1$ experiments are required to settle this matter. In formal terms, one can prove that distinguishability of parameter vectors $x$ from a fixed vector $x_0$ can be attained by means of randomly chosen sets of $r+1$ experiments. The proof of this fact is entirely analogous to the one given for our main theorem; the only difference is that in the definition of the sets $\mathcal{P}_{\beta,q}$, we can simply look at elements $(x, w) \in (\mathbb{X} \setminus \{x_0\}) \times \mathbb{U}^{(q)}$, and this has dimension $qm + r - q$, so using $q = r+1$ we obtain a projection $\mathcal{B}_{\beta,r+1}$, now defined in terms of existence of $x$ not equivalent to $x_0$, and this is a union of manifolds of positive codimension in $\mathbb{U}^{(r+1)}$.